\documentclass[10pt]{amsart}
\usepackage[cp1251]{inputenc}
\usepackage[english,russian]{babel}
\usepackage{amsmath}
\usepackage{amssymb}
\usepackage{amsfonts}
\usepackage{graphicx}

\newtheorem{lemma}{Lemma}
\newtheorem{theorem}{Theorem}

\newtheorem{proposition}{Proposition}

\begin{document}
\renewcommand{\refname}{References}
\renewcommand{\proofname}{Proof.}
\renewcommand{\figurename}{Fig.}

\title[Finite simple groups]{Finite simple groups with two maximal subgroups of coprime orders}
\author{{N.V. Maslova}}%
\address{Natalia Vladimirovna Maslova
\newline\hphantom{iii} Krasovskii Institute of Mathematics and Mechanics UB RAS,
\newline\hphantom{iii} S. Kovalevskaya Str., 16,
\newline\hphantom{iii} 620108, Yekaterinburs, Russia
\newline\hphantom{iii} Ural Federal University,
\newline\hphantom{iii} Turgeneva Str., 4,
\newline\hphantom{iii} 620075, Yekaterinburs, Russia}%
\email{butterson@mail.ru}%

\thanks{\sc Maslova, N.V.,
Finite simple groups with two maximal subgroups of coprime orders}
\thanks{\copyright \ 2022 Maslova N.V.}
\thanks{\rm The reported study was funded by RFBR and BRFBR, project number  20-51-00007 }
%\thanks{\it Received  ..., ..., 2022, published  ..., ...,  202....}%

%\semrtop
\vspace{1cm}
\maketitle

\centerline{\small Dedicated to the memory of Vyacheslav Aleksandrovich Belonogov}

\medskip

{\small
\begin{quote}
\noindent{\sc Abstract. } In 1962, V.~A. Belonogov proved that if a finite group $G$ contains two maximal subgroups of coprime orders, then either $G$ is one of known solvable groups or $G$ is simple. In this short note based on results by M.~Liebeck and J.~Saxl on odd order maximal subgroups in finite simple groups we determine possibilities for triples $(G,H,M)$, where $G$ is a finite nonabelian simple group, $H$ and $M$ are maximal subgroups of $G$ with $(|H|,|M|)=1$.
\medskip

\noindent{\bf Keywords:} finite group, simple group, maximal subgroup, subgroups of coprime orders.
 \end{quote}
}

\section{Introduction}

Throughout the paper we consider only finite groups, and henceforth the term group means
finite group. Our terminology and notation are mostly standard and can be found in \cite{BrHoDou,Atlas,KlLi}.
We denote the group $PSL_n(q)$ by $PSL_n^+(q)$ and the group $PSU_n(q)$ by $PSL_n^-(q)$.
The largest integer power of a prime $p$ dividing a positive integer $k$ is called the
{\it $p$-part} of $k$ and is denoted by $k_p$.
Given an integer $a$ and a positive integer $n$ coprime to $a$, the {\it multiplicative order} $ord_n(a)$ of $a$ modulo $n$ is the smallest positive integer $k$ with $$a^k \equiv 1 \pmod{n}.$$ In other words, the multiplicative order of $a$ modulo $n$ is the order of $a$ in the multiplicative group of the units in the ring of the integers modulo $n$.

Following~\cite{Belonogov}, we write that a group $G$ is a {\it group of type~$A$} if $G$ is a non-special group of order $p q^\beta$, where $p$ and $q$ are primes, a subgroup of order $q^\beta$ is a normal elementary abelian subgroup in $G$, $q^\beta\equiv 1 \pmod{p}$, and $\beta$ is the least with this property.

In 1962, V.~A.~Belonogov~\cite[Theorem~8]{Belonogov} proved the following proposition.

\begin{proposition} If a finite group $G$ contains two maximal subgroups of coprime orders, then $G$ is a cyclic group of order $pq$, where $p$ and $q$ are distinct primes, or $G$ is a group of type~$A$ or $G$ is simple.
\end{proposition}

It is clear that if a simple group $G$ has two maximal subgroups $H$ and $M$ with $(|H|,|M|)=1$, then $|H|$ is odd or $|M|$ is odd. Maximal subgroups of odd orders in simple groups were described by M.~Liebeck and J.~Saxl~\cite{LiSa1}. M.~Aschbacher~\cite{Aschbacher} described natural geometrically defined eight families $\mathfrak{C}_i$ for $1 \le i \le 8$ of subgroups of simple classical groups, which are now called Aschbacher classes, and has proved that if a maximal subgroup of a simple classical group does not belong to the union of Aschbacher classes of the group, then this maximal subgroup is almost simple. Based on the classification of finite simple groups, results by M.~Liebeck and J.~Saxl~\cite{LiSa1}, Aschbacher's classification and further results on maximal subgroups in almost simple groups~\cite{BrHoDou,KlLi} in this short note we prove the following theorem.

\begin{theorem}\label{MThMmAS}
Let $G$ be a finite nonabelian simple group.

$(i)$ If $G$ is sporadic, then $G$ contains two maximal subgroups $H$ and $M$ with $(|H|,|M|)=1$ if and only if one of the following statements holds{\rm:}

$(1)$ $G \cong M_{23}$, $H\cong 23:11$, and $M$ is one of the following subgroups{\rm:} ${PSL_3(4):{2_2}}$, $2^4:A_7$, $A_8$, $2^4:(3\times A_5):2${\rm;}

$(2)$ $G \cong F_2=B$, $H \cong 47:23$ and $M$ is one of the following subgroups{\rm:} ${2^.({^2}E_6(2)):2}$, ${2^{9+16}.PSp_8(2)}$, ${Th}$, ${(2^2 \times F_4(2)):2}$, ${2^{2+10+20}.(M_{22}:2 \times S_3)}$, \\ ${2^{5+5+10+10}.PSL_5(2)}$, ${S_3 \times Fi_{22}:2}$, ${[2^{35}].(S_5 \times PSL_3(2))}$, ${HN:2}$, ${P\Omega_8^+(3):S_4}$, \\ ${3^{1+8}:{2^{1+6}}^.PSU_4(2).2}$, ${5:4 \times HS:2}$, ${(3^2:D_8 \times PSU_4(3).2^2).2}$,  ${S_4 \times {^2}F_4(2)}$, \\ ${3^{2+3+6}.(S_4 \times 2S_4)}$, ${S_5 \times M_{22}:2}$, ${{5^3}^.PSL_3(5)}$, ${(S_6 \times PSL_3(4):2).2}$,  ${5^{1+4}:2^{1+4}.A_5.4}$, ${(S_6 \times S_6).4}$, ${5^2:4S_4 \times S_5}$, ${PSL_2(49)^.2}$, ${M_{11}}$, ${PSL_2(31)}$,  ${PSL_3(3)}$, ${PSL_2(17):2}$, ${PSL_2(11):2}$.

$(ii)$ If $G$ is an alternating group of degree at least $5$ or an exceptional group of Lie type, then $G$ does not contain two maximal subgroups of coprime orders.

$(iii)$ Let $G$ be a simple classical group of Lie type. Assume that $G$ contains two maximal subgroups $H$ and $M$ with $(|H|,|M|)=1$,
and without loss of generality assume that $|H|$ is odd. Then one of the following statements holds{\rm:}

$(1)$ $G \cong PSL_2(q)$, where $7< q \equiv 3 \pmod{4}$, $H \cong E_q:\frac{q-1}{2}$, and $M \cong D_{2(q+1)}${\rm;}

$(2)$ $G \cong PSL_2(q)$, where $q$ is prime and $q \equiv 23 \pmod{24}$, $H \cong E_q:\frac{q-1}{2}$, and $M \cong S_4${\rm;}

$(3)$ $G \cong PSL_2(q)$, where $q$ is prime and $q \equiv 83, 107 \pmod{120}$, $H \cong E_q:\frac{q-1}{2}$, and $M \cong A_4${\rm;}

$(4)$ $G \cong PSL_2(q)$, where $q$ is prime and $q \equiv 59 \pmod{60}$, $H \cong E_q:\frac{q-1}{2}$, and $M \cong A_5${\rm;}

$(5)$ $G \cong PSL_n(q)$, where $n$ is an odd prime, $q$ is not a power of $n$, and $ord_n(q)=n-1$, ${H \cong \frac{q^n-1}{(q-1)(q-1,n)}:n}$, and $M$ is the stabilizer in $G$ of a subspace of dimension $m$ of the natural projective module of $G$, where $2 \le m \le n-2${\rm;}

$(6)$ $G \cong PSU_n(q)$, where $n$ is an odd prime and $q$ is not a power of $n$, \\ ${H \cong \frac{q^n+1}{(q+1)(q+1,n)}:n}$, $M$ is the stabilizer in $G$ of a non-degenerate subspace of dimension $m$ of the natural projective module of $G$, where $2 \le m \le \frac{n-1}{2}$, and one of the following statements holds{\rm:}

$\mbox{ }$$\mbox{ }$$\mbox{ }$$(6a)$  $n \equiv 1 \pmod{4}$, $ord_n(q)=n-1$, and $(n,q)\not = (5,2)${\rm;}

$\mbox{ }$$\mbox{ }$$\mbox{ }$$(6b)$  $n \equiv 3 \pmod{4}$ and $ord_n(q)=\frac{n-1}{2}$.
%$(n,q) \not \in \{(3,5), (5,2)\}$,

$(7)$ $G \cong PSU_3(q)$, where $q\not =5$ and $(q+1)_3=3$, ${H \cong \frac{q^3+1}{(q+1)(q+1,3)}:3}$, and $M$ is the stabilizer in $G$ of a totally singular subspace of dimension $1$ of the natural projective module of $G$.

$(8)$ $G \cong PSU_n(q)$, where $n \ge 5$ is a prime, $q$ is not a power of $n$, $(n,q)\not = (5,2)$, ${H \cong \frac{q^n+1}{(q+1)(q+1,n)}:n}$, $M$ is the stabilizer in $G$ of a totally singular subspace of dimension $m$ of the natural projective module of $G$, and one of the following statements holds{\rm:}

$\mbox{ }$$\mbox{ }$$\mbox{ }$$(8a)$  $n \equiv 1 \pmod 4$, $ord_n(q)=n-1$, and $1 \le m <\frac{n}{3}${\rm;}

$\mbox{ }$$\mbox{ }$$\mbox{ }$$(8b)$  $n \equiv 3 \pmod 4$, $ord_n(q)=\frac{n-1}{2}$, and $1 \le m <\frac{n}{3}${\rm;}

$\mbox{ }$$\mbox{ }$$\mbox{ }$$(8c)$  $ord_n(q)=n-1$ and $\frac{n}{3} < m < \frac{n}{2}${\rm;}

$\mbox{ }$$\mbox{ }$$\mbox{ }$$(8d)$  $n \equiv 3 \pmod{4}$, $ord_n(q)=n-1$, and $\frac{n+1}{4} < m < \frac{n}{3}$.

%where $1 \le m \le \frac{n-3}{2}${\rm;}

$(9)$ $G \cong PSL_n^\varepsilon(q)$, where $n \ge 13$ is a prime and $\varepsilon \in \{+,-\}$, ${H \cong \frac{q^n-\varepsilon1}{(q-\varepsilon1)(q-\varepsilon1,n)}:n}$, and $M$ is almost simple.

Moreover, $H$ is always a subgroup from the union of Aschbacher classes of $G$, and all the possibilities when $M$ is a subgroup from the union of Aschbacher classes of $G$ are listed in Statements $(iii)(1)$--$(iii)(8)$.

\end{theorem}

\section{Proof of Theorem~\ref{MThMmAS}}

The following assertion is well-known and easy-proving.

\begin{lemma}\label{PrPow+-1} Let $q > 1$ be an integer and $k$ and $m$ be positive integers.
Then

$(i)$ $(q^k-1, q^m-1) = q^{(k,m)}-1${\rm;}

$(ii)$ $(q^k+1, q^m+1) = q^{(k,m)}+1$ if $k_2 = m_2$ and $(q^k+1, q^m+1) = (2, q+1)$ otherwise{\rm;}

$(iii)$ $(q^k-1, q^m+1) = q^{(k,m)}+1$ if $k_2 > m_2$ and $(q^k-1, q^m+1) = (2,q+1)$ otherwise.
\end{lemma}

The following result by M.~Liebeck and J.~Saxl will be useful to prove Theorem~\ref{MThMmAS}.

\begin{proposition}[{\rm \cite[Theorem~2]{LiSa1}}]\label{OddOrderMaxSgp} Let $G$ be a finite simple group and $H$ be a maximal subgroup of odd order form $G$.
Then one of the following statements holds{\rm:}

$(1)$ $G \cong A_p$, $H \cong p.\frac{p-1}{2}$, $p$ is prime, $p \equiv 3 \pmod{4}$, and $p \not \in \{7,11,23\}${\rm;}

$(2)$ $G \cong PSL_2(q)$, $H \cong E_q.\frac{q-1}{2}$, and $q \equiv 3 \pmod{4}${\rm;}

$(3)$ $G \cong PSL_n(q)$, $H \cong (\frac{q^n-1}{(q-1)(n,q-1)}).n$, $n$ is an odd prime, and $G \not \cong PSL_3(4)${\rm;}

$(4)$ $G \cong PSU_n(q)$, $H \cong (\frac{q^n+1}{(q+1)(n,q+1)}).n$, $n$ is an odd prime, and $G \not \cong PSU_3(3)$, $PSU_3(5)$, $PSU_5(2)${\rm;}

$(5)$ $G \cong M_{23}$ and $H \cong 23:11${\rm;}

$(6)$ $G \cong Th$ and $H \cong 31:15${\rm;}

$(7)$ $G \cong F_2$ and $H \cong 47:23${\rm;}

$(8)$ $G \cong F_1$ and $H \cong 59:29$ or $H \cong 71:35$.

\end{proposition}

Let $G$ be a nonabelian simple group. Assume that $M$ and $H$ are maximal subgroups of $G$ such that $(|M|,|H|)=1$. Without loss of generality, we can assume that $|H|$ is odd. Now with respect to the classification of finite simple groups, consider simple groups case by case.

\medskip

Let $G$ be sporadic. Then one of Statements $(5)$--$(8)$ of Proposition~\ref{OddOrderMaxSgp} holds.

Assume that $G \cong M_{23}$. Maximal subgroups of $G$ are known, see \cite{Atlas}. Thus, $G$ contains a maximal subgroup $H \cong 23:11$, and $M$ is a maximal subgroup of $G$ such that $(|M|, |H|)=1$ if and only if $M \in \{$ $PSL_3(4):{2_2}$, $2^4:A_7$, $A_8$, $2^4:(3\times A_5):2\}$.

Assume that $G \cong Th$. Maximal subgroups of $G$ are known, see \cite{Eatlas}. $G$ contains a maximal subgroup $H \cong 31:15$ but if $M$ is a maximal subgroup of $G$ which is not $G$-conjugate to $H$, then $|M|$ is divisible by $3$ or $5$, therefore $G$ does not contain a pair of maximal subgroup of coprime orders.

Let $G \cong F_2$ be the Baby Monster group. Maximal subgroups of $G$ are known~\cite{Wilson1999}, their orders are presented also in \cite{Eatlas}. By~\cite{Eatlas,Wilson1999}, $G$ contains a maximal subgroup $H \cong 47:23$, and $M$ is a maximal subgroup of $G$ such that $(|M|, |H|)=1$ if and only if $M$ is one of the following subgroups of $G${\rm:} ${2^.({^2}E_6(2)):2}$, ${2^{9+16}.PSp_8(2)}$, ${Th}$, ${(2^2 \times F_4(2)):2}$, ${2^{2+10+20}.(M_{22}:2 \times S_3)}$, ${2^{5+5+10+10}.PSL_5(2)}$, ${S_3 \times Fi_{22}:2}$, ${[2^{35}].(S_5 \times PSL_3(2))}$, ${HN:2}$, ${P\Omega_8^+(3):S_4}$, ${3^{1+8}:{2^{1+6}}^.PSU_4(2).2}$, ${5:4 \times HS:2}$, ${(3^2:D_8 \times PSU_4(3).2^2).2}$,  ${S_4 \times {^2}F_4(2)}$, ${3^{2+3+6}.(S_4 \times 2S_4)}$, ${S_5 \times M_{22}:2}$, ${{5^3}^.PSL_3(5)}$, ${(S_6 \times PSL_3(4):2).2}$,  ${5^{1+4}:2^{1+4}.A_5.4}$, ${(S_6 \times S_6).4}$, ${5^2:4S_4 \times S_5}$, ${PSL_2(49)^.2}$, ${M_{11}}$, ${PSL_2(31)}$,  ${PSL_3(3)}$, ${PSL_2(17):2}$, ${PSL_2(11):2}$.

Let $G \cong F_1$ be the Monster group and $H_1\cong 59:29$ be a subgroup of $G$. Then $H_1$ is a $\{29,59\}$-Hall subgroup of $G$, and by \cite[Theorem~A]{Gross}, the $\{29,59\}$-Hall subgroups form a unique conjugacy class in $G$. Note that by \cite{Eatlas}, $G$ contains a maximal subgroup $H \cong PSL_2(59)$, and $H$ contains a parabolic maximal subgroup which is a $\{29,59\}$-Hall subgroup of $G$. Thus, $H_1$ is not maximal in $G$ and $G$ does not contain a maximal subgroup of the form $59:29$. Let $H_2\cong 71:35$ be a subgroup of $G$. By~\cite{Wilson1988}, $H_2$ is the normalizer in $G$ of a Sylow $71$-subgroup of $G$, therefore the subgroups of the form $71:35$ form a unique conjugacy class in $G$. Now by \cite{Eatlas}, $G$ contains a maximal subgroup $K \cong PSL_2(71)$, and $K$ contains a parabolic maximal subgroup of the form $71:35$. Thus, $H_2$ is not maximal in $G$ and $G$ does not contain a maximal subgroup of the form $71:35$. 

Statement~$(i)$ of Theorem~\ref{MThMmAS} holds.

\medskip

Let $G$ be a simple alternating group. Then Statement~$(1)$ of Proposition~\ref{OddOrderMaxSgp} holds. Assume that $G$ acts naturally on the
set $\Omega=\{1,\ldots,p\}$. Let $M$ be a maximal subgroup of $G$ such that $p$ does not divide $|M|$. Then $M$ is intransitive on $\Omega$, and  therefore $M$ is the stabilizer in $G$ of a partition $\Omega = \Omega_1 \cup \Omega_2$ of $\Omega$ into disjoint subsets $\Omega_1$ and $\Omega_2$. Now since $p$ is prime, $p \equiv 3 \pmod{4}$, and $p \not \in \{7,11,23\}$, we have that $p \ge 19$, and it is clear that $|M|$ is not coprime to $\frac{p-1}{2}$. Therefore $G$ does not contain a pair of maximal subgroup of coprime orders. Thus, Statement~$(ii)$ of Theorem~\ref{MThMmAS} holds for alternating groups.

\medskip

If $G$ is an exceptional group of Lie type, then by Proposition~\ref{OddOrderMaxSgp}, $G$ does not contain a maximal subgroup of odd order. Thus, Statement~$(ii)$ of Theorem~\ref{MThMmAS} holds for exceptional groups of Lie type.

\medskip

Let $G$ be a finite simple classical group. Then one of Statements~$(2)$--$(4)$ of Proposition~\ref{OddOrderMaxSgp} holds.

Let Statement~$(2)$ of Proposition~\ref{OddOrderMaxSgp} hold, i.\,e. $G=PSL_2(q)$ for $3<q \equiv 3 \pmod{4}$. Maximal subgroups of $G$ are known (see, for example, \cite[Tables~8.1,~8.2,~8.7]{BrHoDou}). Note that the group $G=PSL_2(q)$ contains a parabolic maximal subgroup $H$ of the form $E_q.\frac{q-1}{2}$, and $|H|$ is odd if and only if $q \equiv 3 \pmod{4}$. We consider possibilities for $M$ case by case.

If $M \cong D_{2(q-1)}$, then it is clear that $(|M|,|H|)\not =1$.

Let $M \cong D_{2(q+1)}$. Then $(|H|,|M|)=1$ if and only if $q \equiv 3 \pmod{4}$ and $M$ is maximal in $G$ if and only if $q \not \in \{7,9\}$.
Thus, Statement~$(iii)(1)$ of Theorem~\ref{MThMmAS} holds.

Let $M \cong S_4$.
Note that $M$ is maximal in $G$ if and only if $q$ is prime and $q \equiv \pm 1 \pmod {8}$, in particular, $q$ is not a power of $3$.
Then $(|M|,|H|)=1$ if and only if $q \equiv 3 \pmod{4}$ and $q\equiv 2\pmod{3}$, i.\,e. $q \equiv 11 \pmod{12}$.
Therefore $M$ is a maximal subgroup of $G$ such that $(|M|,|H|)=1$ if and only if $q$ is a prime and $q \equiv 23 \pmod{24}$.
%11, 23, 35, 47, 59,  71, 83, 95, 107, 119
%11  -17  -5  7  19  -19  3   15  -13   -1
Thus, Statement~$(iii)(2)$ of Theorem~\ref{MThMmAS} holds.

Let $M \cong A_4$. Note that $M$ is maximal in $G$ if and only if either $q=5$ or $q$ is prime and $q \equiv \pm 3, \pm 13 \pmod {40}$. Then again $(|M|,|H|)=1$ if and only if $q \equiv 11 \pmod{12}$.
 Therefore $M$ is a maximal subgroup of $G$ such that $(|M|,|H|)=1$ if and only if $q$ is a prime and $q \equiv 83, 107 \pmod{120}$.
%11, 23, 35, 47, 59,  71, 83, 95, 107, 119
%11  -17  -5  7  19  -19  3   15  -13   -1
Thus, Statement~$(iii)(3)$ of Theorem~\ref{MThMmAS} holds.

If $M \cong PSL_2(q_0).2$, where $q=q_0^2$, or $M \cong PSL_2(q_0)$, where $q=q_0^r$ and $r$ is an odd prime, then it is clear that $(|M|,|H|)\not =1$.

Let $M \cong A_5$. Note that $G=PSL_2(q)$ contains a maximal subgroup isomorphic to $A_5$ if and only if either $q$ is prime and $q\equiv \pm 1 \pmod{10}$ or $q=p^2$, where $p$ is prime and $p \equiv \pm 3 \pmod{10}$. It is easy to see that $(|M|,|H|)=1$ if and only if $q \equiv 11 \pmod{12}$ and $q \equiv 2,3,4 \pmod{5}$. Note that if $q=p^2$, where $p \equiv \pm 3 \pmod{10}$, then $q \equiv 9, 49 \pmod{60}$, therefore this case does not appear. So, $M$ is a maximal subgroup of $G$ such that $(|M|,|H|)=1$ if and only if $q$ is prime and $q \equiv 59 \pmod{60}$.
%q=p \equiv \pm 1 \pmod 10 => q\equiv -1 \pmod 10 and q\equiv -1 \pmod 11 => q \equiv 59 \pmod 60
% mod 60:     9  19  29  39  49  59
% mod 12 =>   9   7   5   3   1  11
%
% q=p^2, p \equiv \pm 3 \pmod 10
% p \pmod 60: 3,    7,   13,   17,   23,   27,   33,   37,   43,   47,   53,   57
% p^2 mod 60: 9,   49,   49,   49,   49,    9,    9,   49,   49,   49,   49,    9
% p^2 mod 12: 9,    1,    1,    1,    1,    9,    9,    1,    1,    1,    1,    9
Thus, Statement~$(iii)(4)$ of Theorem~\ref{MThMmAS} holds.

\medskip

Let Statement~$(3)$ of Proposition~\ref{OddOrderMaxSgp} hold, i.\,e. $G \cong PSL_n(q)$, $H \cong (\frac{q^n-1}{(q-1)(n,q-1)}).n$, $n$ is an odd prime, and $G \not \cong PSL_3(4)$. Note that by \cite[Table~8.3]{BrHoDou}, the group $PSL_3(4)$ does not contain a maximal subgroup of odd order. Thus, we can assume that $G \not \cong PSL_3(4)$, and by \cite[Tables~8.3,~8.18,~8.35,~8.70]{BrHoDou} and \cite[Table~3.5.A, Proposition~4.3.6]{KlLi}, $G$ has a maximal subgroup $H \cong (\frac{q^n-1}{(q-1)(n,q-1)}).n$ such that $|H|$ is odd. By the Aschbacher theorem, for each maximal subgroup $M$ of $G$, if $M$ does not belong to the union of Aschbacher classes $\mathfrak{C}_i(G)$ for $i \in \{1,\ldots, 8\}$, then $M$ is almost simple and we write $M \in \mathfrak{S}(G)$ in this case. If $n \le 11$, then all maximal subgroups of $G$ are known, see \cite[Tables~8.3,~8.4,~8.18,~8.19,~8.35,~8.36,~8.70~8.71]{BrHoDou}, and if $n \ge 13$, then maximal subgroups of $G$ from the union of Aschbacher classes of $G$ are known, see \cite[Table~3.5.A]{KlLi}. So, we consider possibilities for $M$ case by case.

Since $n$ is odd and prime, by Tables~8.3,~8.4,~8.18,~8.19,~8.35,~8.36,~8.70, and~8.71 of \cite{BrHoDou} and \cite[Table~3.5.A, Propositions~4.1.17,~4.2.9,~4.5.3,~4.6.5,~4.8.4,~4.8.5]{KlLi}, if $M$ is a maximal subgroup of $G$ which is not isomorphic to $H$, then one of the following statements holds{\rm:}

\begin{itemize}
\item[$(\mathfrak{C}_1)$] $M$ is the stabilizer in $G$ of a subspace of dimension $m$ of the natural projective module of $G$, where $1 \le m \le n-1${\rm;}
\item[$(\mathfrak{C}_2)$] $M \cong [\frac{(q-1)^{n-1}}{(q-1,n)}]:S_n$ is the stabilizer in $G$ of a decomposition of the natural projective module of $G$ into a direct sum of subspaces of dimension $1${\rm;}
\item[$(\mathfrak{C}_5)$] $M \cong \frac{c}{(q-1,n)}PGL_n(q_0)$ if $q=q_0^r$ for odd prime $r$, where $c=\frac{(q-1)}{[q_0-1,\frac{q-1}{(q-1,n)}]}${\rm;}
\item[$(\mathfrak{C}_6)$] $M \cong [n^3].Sp_2(n)$ if some special conditions for $n$ and $q$ hold and $3^2.Q_8$ if $n=3$ and some special conditions for $q$ hold{\rm;}
\item[$(\mathfrak{C}_8)$] $M \cong PSO_n(q)$ or $M \cong PSU_n(q_0).[\frac{(q_0+1,n)c}{(q-1,n)}]$ for $c=\frac{(q-1)}{[q_0+1,\frac{q-1}{(q-1,n)}]}$ if $q=q_0^2${\rm;}
\item[$(\mathfrak{S})$] $M$ is almost simple.
\end{itemize}

Note that by~\cite[Tables~8.3,~8.18,~8.35,~8.70]{BrHoDou} and~\cite[Proposition~4.1.17]{KlLi}, for each $1 \le m \le n-1$, the stabilizer $M$ in $G$ of a subspace of dimension $m$ of the natural projective module of $G$ is a maximal subgroup of $G$. Moreover, by \cite[Proposition~4.1.17]{KlLi}, $$|M|=\frac{q^{n(n-1)}\cdot (q-1)}{(q-1,n)}\cdot\prod_{i=2}^m(q^i-1)\cdot \prod_{i=2}^{n-m}(q^i-1).$$ By Lemma~\ref{PrPow+-1}, is clear that $(|M|,|H|)=1$ if and only if $n$ does not divide $|M|$. Note that by Fermat’s little theorem, if $q$ is not a power of $n$, then $n$ divides $q^{n-1}-1$. So, $n$ does not divide $|M|$ if and only if $q$ is not a power of $n$, $2 \le m \le n-2$, and $ord_n(q)=n-1$. Thus, Statement~$(iii)(5)$ of Theorem~\ref{MThMmAS} holds.

If $M \cong [\frac{(q-1)^{n-1}}{(q-1,n)}]:S_n$, then it is clear that $n$ divides $|M|$, therefore $|M|$ and $|H|$ are not coprime.

If  $M \cong \frac{c}{(q-1,n)}PGL_n(q_0)$ for $q=q_0^r$ with $r$ odd prime and $c=\frac{(q-1)}{[q_0-1,\frac{q-1}{(q-1,n)}]}$ or $M \cong PSU_n(q_0).[\frac{(q_0+1,n)c}{(q-1,n)}]$ for $q=q_0^2$ and $c=\frac{(q-1)}{[q_0+1,\frac{q-1}{(q-1,n)}]}$, then $q_0^{n-1}-1$ divides $|M|$, and by Fermat’s little theorem, if $q$ is not a power of $n$, then $n$ divides $q_0^{n-1}-1$. Thus, in any case $n$ divides $|M|$, therefore $|M|$ and $|H|$ are not coprime.

If $M \cong [n^3].Sp_2(n)$, then it is clear that $n$ divides $|M|$, therefore $|M|$ and $|H|$ are not coprime.

If $M \cong PSO_n(q)$, then $q^{n-1}-1$ divides $|M|$, again using Fermat’s little theorem, we conclude that $n$ divides $|M|$. Thus, $|M|$ and $|H|$ are not coprime.

Almost simple subgroups of low-dimensional finite simple linear groups are known, see \cite[Tables~8.4,~8.19,~8.36,~8.71]{BrHoDou}. Following \cite{BrHoDou}, here we provide the preimage $\tilde{M}$ in $\tilde{G}=SL_n(q)$ if in some cases $G$ has a maximal almost simple subgroup isomorphic to $M$.
\medskip

\centerline{ ''Almost simple'' maximal subgroups of $SL_n(q)$ for $n \in \{3,5,7,11\}$}

\smallskip

\centerline{\begin{tabular}{|p{10mm}|p{20mm}|p{100mm}|}
\hline
$n$ & $\mbox{Group } \tilde{G}$ & $\mbox{ }$$\mbox{ }$$\mbox{ }$$\mbox{ }$$\mbox{ }$$\mbox{ }$$\mbox{ }$$\mbox{ }$$\mbox{ }$$\mbox{ }$$\mbox{ }$$\mbox{ }$$\mbox{ }$$\mbox{ }$$\mbox{ }$$\mbox{ }$$\mbox{ }$$\mbox{ }$$\mbox{ }$ $\mbox{Subgroups } \tilde{M}$ \\ \hline

$3$ &$SL_3(q)$& $(q-1,3) \times PSL_2(7)$, $3^.A_6 $\\ \hline

$5$ &$SL_5(q)$& $(q-1,5) \times PSL_2(11)$, $M_{11}$, $(q-1,5) \times PSU_4(2)$\\ \hline

$7$ &$SL_7(q)$& $(q-1,7) \times PSU_3(3)$\\ \hline

$11$ &$SL_{11}(q)$& ${PSL_2(23)}$, ${M_{24}}$, ${(q-1,11) \times PSL_2(23)}$, ${(q-1,11) \times PSU_5(2)}$\\ \hline

\end{tabular}

}

\medskip

Now it is clear that if $G=PSL_n(q)$ for $n \le 11$ and $M$ is a maximal almost simple subgroup from $\mathfrak{S}(G)$, then $n$ divides $|M|$, therefore $|M|$ and $|H|$ are not coprime. Thus, Statement~$(iii)(9)$ of Theorem~\ref{MThMmAS} holds for simple linear groups.

\medskip

Let Statement~$(4)$ of Proposition~\ref{OddOrderMaxSgp} hold, i.\,e. $G \cong PSU_n(q)$, $H \cong (\frac{q^n+1}{(q+1)(n,q+1)}).n$, $n$ is an odd prime, and $G \not \cong PSU_3(3)$, $PSU_3(5)$, $PSU_5(2)$. Note that by \cite[Tables~8.5,~8.20]{BrHoDou} groups $PSU_3(3)$, $PSU_3(5)$, and $PSU_5(2)$ do not contain maximal subgroups of odd order. Thus, we can assume that $G \not \cong PSU_3(3)$, $PSU_3(5)$, $PSU_5(2)$, and by \cite[Tables~8.5,~8.20,~8.37,~8.72]{BrHoDou} and \cite[Table~3.5.B, Proposition~4.3.6]{KlLi}, $G$ has a maximal subgroup $H \cong (\frac{q^n+1}{(q+1)(n,q+1)}).n$ such that $|H|$ is odd. By the Aschbacher theorem, for each maximal subgroup $M$ of $G$, if $M$ does not belong to the union of Aschbacher classes $\mathfrak{C}_i(G)$ for $i \in \{1,\ldots, 8\}$, then $M$ is almost simple and we again write $M \in \mathfrak{S}(G)$ in this case. If $n \le 11$, then all maximal subgroups of $G$ are known, see \cite[Tables~8.5,~8.6,~8.20,~8.21, 8.37, 8.38, 8.72,~8.73]{BrHoDou}, and if $n \ge 13$, then maximal subgroups of $G$ from the union of Aschbacher classes of $G$ are known, see \cite[Table~3.5.B]{KlLi}. So, we consider possibilities for $M$ case by case.

Since $n$ is odd and prime, by Tables~8.5,~8.6,~8.20,~8.21,~8.37,~8.38,~8.72,~8.73 of \cite{BrHoDou} and \cite[Table~3.5.B, Propositions~4.1.4,~4.1.18,~4.2.9,~4.5.3,~4.5.5,~4.6.5]{KlLi}, if $M$ is a maximal subgroup of $G$ which is not isomorphic to $H$, then one of the following statements holds{\rm:}

\begin{itemize}
\item[$(\mathfrak{C}_1)$] $M$ is the stabilizer in $G$ of a subspace of dimension $m$, where $1 \le m \le \frac{n-1}{2}$, of the natural projective module of $G$, and either the subspace is totally singular or the subspace is non-degenerate{\rm;}
\item[$(\mathfrak{C}_2)$] $M \cong [\frac{(q+1)^{n-1}}{(q+1,n)}]:S_n$ is the stabilizer in $G$ of a decomposition of the natural projective module of $G$ into a direct sum of pairwise orthogonal non-degenerate subspaces of dimension $1${\rm;}
\item[$(\mathfrak{C}_5)$] $M \cong \frac{c}{(q+1,n)}PGU_n(q_0)$ if $q=q_0^r$ for odd prime $r$, where $c=\frac{(q+1)}{[q_0+1,\frac{q+1}{(q+1,n)}]}$, or $M \cong PSO_n(q)${\rm;}
\item[$(\mathfrak{C}_6)$] $M \cong [n^3].Sp_2(n)$ if some special conditions for $n$ and $q$ hold or $3^2.Q_8$ if $n=3$ and some special conditions for $q$ hold{\rm;}
\item[$(\mathfrak{S})$] $M$ is almost simple.
\end{itemize}

Note that by~\cite[Tables~8.5,~8.20,~8.37,~8.72]{BrHoDou} and~\cite[Proposition~4.1.4]{KlLi}, for each $1 \le m \le \frac{n-1}{2}$, the stabilizer $M$ in $G$ of a non-degenerate subspace of dimension $m$ of the natural projective module of $G$ is a maximal subgroup of $G$. Moreover, by \cite[Proposition~4.1.4]{KlLi}, $$|M|=\frac{q^{\frac{m^2+(n-m)^2-n}{2}}\cdot (q+1)}{(q+1,n)}\cdot\prod_{i=2}^m(q^i-(-1)^i)\cdot \prod_{i=2}^{n-m}(q^i-(-1)^i).$$ By Lemma~\ref{PrPow+-1}, is clear that $(|M|,|H|)=1$ if and only if $n$ does not divide $|M|$. Note that by Fermat’s little theorem, if $q$ is not a power of $n$, then $n$ divides $q^{n-1}-1$, therefore $ord_n(q)$ divides $n-1$ and if $n$ does not divide $|M|$, then $m \ge 2$. Let $n-1=k\cdot ord_n(q)$. If $k \ge 3$ is odd, then $n$ divides $q^{\frac{n-1}{k}}-1$ and $q^{\frac{n-1}{k}}-1$ divides $|M|$ since $\frac{n-1}{k}<\frac{n-1}{2}$ and $\frac{n-1}{k}$ is even. If $k \ge 4$ is even, then $n$ divides $q^{\frac{2 \cdot (n-1)}{k}}-1$ and $q^{\frac{2 \cdot (n-1)}{k}}-1$ divides $|M|$ since $\frac{2(n-1)}{k} \le \frac{n-1}{2}$ and $\frac{2(n-1)}{k}$ is even. Thus, if $n$ does not divide $|M|$, then $k \in \{1, 2\}$.

Assume that $n \equiv 1 \pmod 4$, therefore $\frac{n-1}{2}$ is even and $q^\frac{n-1}{2}-1$ divides $|M|$. Thus, in this case if $k=2$, then $n$ divides $|M|$. Suppose that $k=1$ and $n$ divides $|M|$. Then $n$ divides $q^i+1$ for some odd $i \le n-2$. Therefore $n-1=ord_n(q)$ divides $2i \le 2n-4$. Thus, $2i=n-1$ and $i$ is even, a contradiction.

Assume that $n \equiv 3 \pmod 4$, therefore $\frac{n-1}{2}$ is odd and $q^\frac{n-1}{2}+1$ divides $|M|$. Thus, in this case if $k=1$, then $n$ divides $|M|$. Suppose that $k=2$ and $n$ divides $|M|$. Then since $\frac{n-1}{2}$ is odd, we have $q^\frac{n-1}{2}-1$ does not divide $|M|$, therefore $n$ divides $q^i+1$ for some odd $i \le n-2$. Thus, $\frac{n-1}{2}=ord_n(q)$ divides $2i \le 2n-4$. We have $2i=n-1$ and $i=\frac{n-1}{2}$, a contradiction. Thus, Statement~$(iii)(6)$ of Theorem~\ref{MThMmAS} holds.

Note that by~\cite[Tables~8.5,~8.20,~8.37,~8.72]{BrHoDou} and~\cite[Proposition~4.1.18]{KlLi}, for each $1 \le m \le \frac{n-1}{2}$, the stabilizer $M$ in $G$ of a totally singular subspace of dimension $m$ of the natural projective module of $G$ is a maximal subgroup of $G$. Moreover, by \cite[Proposition~4.1.18]{KlLi}, $$|M|=\frac{q^{\frac{n(n-1)}{2}}\cdot (q^2-1)}{(q+1,n)}\cdot\prod_{i=2}^m(q^{2i}-1)\cdot \prod_{i=2}^{n-2m}(q^i-(-1)^i).$$
By Lemma~\ref{PrPow+-1}, is clear that $(|M|,|H|)=1$ if and only if $n$ does not divide $|M|$. Note that by Fermat’s little theorem, if $q$ is not a power of $n$, then $n$ divides $q^{n-1}-1$, therefore $ord_n(q)$ divides $n-1$.

Let $n=3$. Then $m=1$, $|M| = \frac{q^3\cdot (q^2-1)}{(q+1,3)}$, and $3$ does not divide $|M|$ if and only if $(q+1)_3=3$. Thus, Statement~$(iii)(7)$ of Theorem~\ref{MThMmAS} holds.

Let $n\ge 5$. Assume that $m > \frac{n}{3}$. Then the number $\prod_{i=2}^{n-2m}(q^i-(-1)^i)$ divides the number $\prod_{i=2}^m(q^{2i}-1)$, therefore $n$ does not divide $|M|$ if and only if $ord_n(q^2)>m >\frac{n}{3}$. It is clear that $ord_n(q^2)$ divides $\frac{n-1}{2}$, therefore $n$ does not divide $|M|$ if and only if $ord_n(q^2)=\frac{n-1}{2}$ which is equivalent to $ord_n(q)=n-1$.

Assume that $m < \frac{n}{4}$. Then the number $\prod_{i=2}^m(q^{2i}-1)$ divides the number $\prod_{i=2}^{n-2m}(q^i-(-1)^i)$. Thus, as in the proof of Statement~$(iii)(6)$, we obtain that $n$ does not divide $|M|$ if and only if either $n \equiv 1 \pmod{4}$ and $ord_n(q) = n - 1$ or $n \equiv 3 \pmod{4}$ and $ord_n(q) = \frac{n-1}{2}$.

Finally, assume that $\frac{n}{4} < m <\frac{n}{3}$. Let $n-1=k\cdot ord_n(q)$. If $k \ge 3$ and $\frac{n-1}{k}$ is even, then $n$ divides $q^{\frac{n-1}{k}}-1$ and $q^{\frac{n-1}{k}}-1$ divides $\prod_{i=2}^m(q^{2i}-1)$ since $\frac{n-1}{k}<\frac{n-1}{2}<2\cdot \frac{n}{4}$. If $k \ge 4$ and $\frac{n-1}{k}$ is odd, then $n$ divides $q^{\frac{n-1}{k}}-1$ and $q^{\frac{n-1}{k}}-1$ divides $\prod_{i=2}^m(q^{2i}-1)$ since $2\cdot \frac{n-1}{k} \le 2\cdot \frac{n-1}{4}<2\cdot \frac{n}{4}$. Thus, $k \in \{1, 2\}$.

Assume that $n \equiv 1 \pmod 4$, therefore $\frac{n-1}{2}$ is even. Since $\frac{n}{4} < m$, we have $q^\frac{n-1}{2}-1$ divides $|M|$. Thus, in this case if $k=2$, then $n$ divides $|M|$. Suppose that $k=1$ and $n$ divides $|M|$. Then $n$ divides $q^i-1$ for some $i \le \frac{2n}{3}$ or $n$ divides $q^i+1$ for some odd $i \le n-2$. Note that in the latter case $n-1=ord_n(q)$ divides $2i \le 2n-4$, therefore $2i=n-1$ and $i$ is even. In both cases we have a contradiction.

Assume that $n \equiv 3 \pmod 4$, therefore $\frac{n-1}{2}$ is odd. Suppose that $k=2$ and $n$ divides $|M|$. Then since $\frac{n-1}{2}$ is odd, we have $q^\frac{n-1}{2}-1$ does not divide $|M|$, therefore $n$ divides $q^i-1$ for some even $i \le \frac{2n}{3}$ or $n$ divides $q^i+1$ for some odd $i \le n-2$. In the latter case, $\frac{n-1}{2}=ord_n(q)$ divides $2i \le 2n-4$. Thus, $2i=n-1$ and $i=\frac{n-1}{2}$, a contradiction. Assume that $k=1$. If $m=\frac{n+1}{4}$, then $n$ divides $q^\frac{n-1}{2}+1$ and $q^\frac{n-1}{2}+1$ divides $|M|$. Let $m>\frac{n+1}{4}$. Suppose that $n$ divides $|M|$. Then $n$ divides $q^{2i}-1$ for some $i<\frac{n}{3}$ or $n$ divides $q^i-1$ for some even $i<\frac{n-1}{2}$ or $n$ divides $q^i+1$ for some odd $i <\frac{n-1}{2}$. In the latter case, $n$ divides $q^{2i}-1$ and $2i<n-1$. In all cases we have contradictions. Thus, Statement~$(iii)(8)$ of Theorem~\ref{MThMmAS} holds.

If $M \cong [\frac{(q+1)^{n-1}}{(q+1,n)}]:S_n$, then it is clear that $n$ divides $|M|$, therefore $|M|$ and $|H|$ are not coprime.

If  $M \cong \frac{c}{(q+1,n)}PGU_n(q_0)$ for $q=q_0^r$ with $r$ odd prime and $c=\frac{(q+1)}{[q_0+1,\frac{q+1}{(q+1,n)}]}$, then $q_0^{n-1}-1$ divides $|M|$, and by Fermat’s little theorem, if $q$ is not a power of $n$, then $n$ divides $q_0^{n-1}-1$. Thus, in any case $n$ divides $|M|$, therefore $|M|$ and $|H|$ are not coprime.

Similarly, if $M \cong PSO_n(q)$, then $q^{n-1}-1$ divides $|M|$, again using Fermat’s little theorem, we conclude that $n$ divides $|M|$. Thus, $|M|$ and $|H|$ are not coprime.

If $M \cong [n^3].Sp_2(n)$, then it is clear that $n$ divides $|M|$, therefore $|M|$ and $|H|$ are not coprime.

Almost simple subgroups of low-dimensional finite simple unitary groups are known, see \cite[Tables~8.6,~8.21,~8.38,~8.73]{BrHoDou}.
Following \cite{BrHoDou}, here we provide the preimage $\tilde{M}$ in $\tilde{G}=SU_n(q)$ if in some cases $G$ has a maximal almost simple subgroup isomorphic to $M$.

\medskip

\centerline{ ''Almost simple'' maximal subgroups of $SU_n(q)$ for $n \in \{3,5,7,11\}$}

\smallskip

\centerline{\begin{tabular}{|p{10mm}|p{20mm}|p{100mm}|}
\hline
$n$ & $\mbox{Group } \tilde{G}$ & $\mbox{ }$$\mbox{ }$$\mbox{ }$$\mbox{ }$$\mbox{ }$$\mbox{ }$$\mbox{ }$$\mbox{ }$$\mbox{ }$$\mbox{ }$$\mbox{ }$$\mbox{ }$$\mbox{ }$$\mbox{ }$$\mbox{ }$$\mbox{ }$$\mbox{ }$$\mbox{ }$$\mbox{ }$ $\mbox{Subgroups } \tilde{M}$  \\ \hline

$3$ &$SU_3(q)$& $(q+1,3) \times PSL_2(7)$, $3^.A_6 $, $3^.A_6^.2_3$, $3^.A_7$\\ \hline

$5$ &$SU_5(q)$& $(q+1,5) \times PSL_2(11)$, $(q+1,5) \times PSU_4(2)$\\ \hline

$7$ &$SU_7(q)$& $(q+1,7) \times PSU_3(3)$\\ \hline

$11$ &$SU_{11}(q)$& ${(q+1,11) \times PSL_2(23)}$, ${(q+1,11) \times PSU_5(2)}$\\ \hline

\end{tabular}

}

\medskip

Now it is clear that if $G=PSU_n(q)$ for $n \le 11$ and $M$ is a maximal almost simple subgroup from $\mathfrak{S}(G)$, then $n$ divides $|M|$, therefore $|M|$ and $|H|$ are not coprime. Thus, Statement~$(iii)(9)$ of Theorem~\ref{MThMmAS} holds for simple unitary groups.

\medskip

Note that $H$ is always a subgroup from the union of Aschbacher classes of $G$, and all the possibilities when $M$ is a subgroup from the union of Aschbacher classes of $G$ are listed in Statements~$(iii)(1)$--$(iii)(8)$ of Theorem~\ref{MThMmAS}.

\medskip

%\begin{proof}
%Text of the proof.
%\end{proof}

\section{Acknowledgment}
The author is thankful to Prof. Viktor~S. Monakhov who has drawn the author's attention to the problem considered in this short note, to Dr. Alexey~M. Staroletov and to anonymous reviewers for their helpful comments which improved this text.

\smallskip

The author dedicates this work to the memory of Prof. Vyacheslav Aleksandrovich Belonogov who passed away on December 21, 2021. Prof. Belonogov was alert and actively worked up to the last his hours, and the author has had a chance to discuss the results of this work with him.


\bigskip
\begin{thebibliography}{1}

\bibitem{Aschbacher}
M. Aschbacher, {\it On the maximal subgroups of the finite classical groups}, Invent. Math. {\bf 76}:3
(1984), 469--514.


\bibitem{Belonogov}
V.A. Belonogov, {\it On maximal subgroups. II}, Izv. Vyssh. Uchebn. Zaved. Mat., {\bf 5} (1962), 3--11 (in Russian).


\bibitem{BrHoDou}
J.N. Bray, D.F. Holt, C.M. Roney-Dougal, {\it The maximal subgroups of the low-dimensional finite classical groups}, London Math. Soc. Lect. Note Ser., 407, CUP, Cambridge, 2013.

\bibitem{Atlas} J.H. Conway, et al., {\it Atlas of Finite Groups}, Clarendon Press, Oxford, 1985.

\bibitem{Gross}
F. Gross, {\it Conjugacy of odd order Hall subgroups}, Bull. London Math. Soc., \textbf{19}:4 (1987),
311--319.


\bibitem{KlLi}
P.B. Kleidman, M. Liebeck, {\it The subgroup structure of the finite classical groups}, CUP,
Cambridge, 1990.


%\bibitem{LiZhang} C.~H.~Li, H.~Zhang, {\it The finite primitive groups with soluble stabilizers, and the edge-primitive s-arc transitive graphs}, Proc. London Math. Soc. (3), {\bf 103} (2011), 441--472.

\bibitem{LiSa1}
M.~Liebeck, J.~Saxl, {\it On point stabilizers in primitive permutation groups}, Communicatons in Algebra, {\bf 19}:10 (1991), 2777--2786.

%\bibitem{Schroder} A.~K.~Schroder, The maximal subgroups of the classical groups in dimension 13, 14 and 15, Ph. D. Thesis, University of St Andrews, St Andrews, 2015.

\bibitem{Eatlas}
R.~A.~Wilson [et. al.], {\it Atlas of finite group representations}, https://brauer.maths.qmul.ac.uk/Atlas/v3/.


\bibitem{Wilson1988}
R.~A.~Wilson, {\it The local subgroups of the Monster}, J. Austral. Math. Soc. (A), \textbf{44} (1988), 1--16.

\bibitem{Wilson1999}
R.~A.~Wilson, {\it The maximal subgroups of the Baby Monster, I}, J. Algebra, \textbf{211} (1999), 1--14.

%\bibitem{Wilson2009} R.~A.~Wilson, {\it The finite smple groups},  Springer-Verlag, London, 2009.


%\bibitem{Wilson2017} R.~A.~Wilson, {\it Maximal subgroups of sporadic groups}, arXiv:1701.02095 [math.GR].

%\bibitem{1} D.G. Northcott, {\it Ideal Theory}, Cambridge University Press, Cambridge, 1953.

%\bibitem{2} I.I. Ivanov, {\it Some properties of commutative algebras}, Siberian Electronic Mathematical Reports, {\bf 1} (2004), 142--143.

%\bibitem{3} A.V. Tetenov, M. Samuel, D.A. Vauilin, {\it On dendrites generated by polyhedral systems and their ramification points}, Trudy Inst. Mat. Mekh. UrO RAN, \textbf{23}:4 (2017), 281--291.

\end{thebibliography}
\end{document}